\documentclass[a4paper,12pt,reqno]{amsart}
\usepackage{amsmath,amssymb}

\usepackage{graphicx}
\usepackage{amsthm}
\usepackage{thmtools}
\usepackage{comment}
\usepackage{hyperref}
\usepackage{float}
\usepackage{caption}
\usepackage{mathtools}
\usepackage{arcs}
\usepackage{yhmath}
\usepackage{tikz}
\usetikzlibrary{calc,intersections}
\makeatletter

\@namedef{subjclassname@2020}{%
\textup{2020} Mathematics Subject Classification} 
\makeatother

\newtheorem{theorem}{Theorem}

\newtheorem{proposition}[theorem]{Proposition}

\title{A Note on the Converse Sendov Problem}

\author{Dragomir Grozev}
\address{Institute of Mathematics and Informatics, Bulgarian Academy of Sciences, Acad. G. Bonchev 8,
1113 Sofia, Bulgaria}
\email{drago.grozev@gmail.com}

\author{Nikolai Nikolov}
\address{Institute of Mathematics and Informatics, Bulgarian Academy of Sciences, Acad. G. Bonchev 8,
1113 Sofia, Bulgaria
\vspace{1mm}
\newline Faculty of Information Sciences, State University of Library Studies and Information
Technologies, Shipchenski prohod 69A, 1574 Sofia, Bulgaria}
\email{nik@math.bas.bg}

\thanks{The second named author was partially supported by the Bulgarian National Science Fund,
Ministry of Education and Science of Bulgaria under contract
KP-06-N82/6.} 

\begin{document}

\begin{abstract}
For a polynomial of degree $n$ whose zeros lie in the closed unit disk, we
determine the largest possible distance from a prescribed critical point of
modulus $r$ to the nearest zero. If $n$ is even, the sharp radius is
$\sqrt{1-r^2}$; if $n$ is odd, the sharp radius is strictly smaller for
$r\in(0,1)$ and depends on $n$. Equality cases are also determined.
The proof is based on the logarithmic-derivative identity and elementary
geometric considerations.
\end{abstract}

\maketitle
\vspace{-8mm}

\section{Introduction and results}

Sendov's conjecture asserts that if all zeros of a polynomial lie in the
closed unit disk $\mathbb{D}$, then every zero has a critical point within
distance $1$. It was recently proved for all degrees $n\geq2$ by Mazur
\cite{Mazur}; see Tao \cite{Tao} for an exposition. We consider the converse Sendov problem: how close must a zero lie to a prescribed critical point?

For $n\geq2$ and $0\leq r\leq1$, define
\[
R_n(r):=
\sup\left\{
\min_{1\leq j\leq n}|z_j-w|:
p(z)=c\prod_{j=1}^n(z-z_j),\
|z_j|\leq1,\ p'(w)=0,\ |w|=r
\right\}.
\]

Sofi and Shah \cite{SofiShah} proved that for every critical point $w$ there
is a zero $z_j$ such that
\[
|z_j-w|^2\leq |z_j|^2-|w|^2.
\]
Thus, if all zeros lie in $\mathbb D$,
\[
|z_j-w|\leq\sqrt{1-|w|^2}.
\]
We determine the sharp bound when the degree $n$ is fixed, revealing a distinction between even and odd degrees.

\begin{theorem}
\label{thm:thm_1}
Let $p$ be a monic polynomial of degree $n\geq 2$ whose zero set $Z(p)$ lies
in $\mathbb{D}$, and let $w$ be a critical point of $p$, with
$r=|w|$, $0\leq r\leq 1$. If $n$ is even, then
\begin{equation}
\label{eq:eq_01}
\operatorname{dist}(w,Z(p))
\leq E(r):=\sqrt{1-r^2}.
\end{equation}
If $n$ is odd, then
\begin{equation}
\label{eq:eq_02}
\operatorname{dist}(w,Z(p))
\leq O_n(r)
:=
(1+r)
\sqrt{
\frac{(n-1)(1-r)}
     {(n-1)+(n+1)r}
}.
\end{equation}

Both estimates are sharp. More precisely, let $0<r<1$, $w=r$, and
$m=\lfloor n/2\rfloor$. Equality in \eqref{eq:eq_01} occurs precisely for
\[
p(z)
=
\left(z-r-i\sqrt{1-r^2}\right)^m
\left(z-r+i\sqrt{1-r^2}\right)^m,
\]
whereas equality in~\eqref{eq:eq_02} occurs precisely for
\[
p(z)
=
(z+1)(z-e^{i\theta})^m(z-e^{-i\theta})^m,
\]
where
\[
\cos\theta
=
\frac{(2m+1)r^2+2mr+1}
     {2\bigl(m+(m+1)r\bigr)}.
\]

For $r=0$, equality holds precisely when all zeros lie on
$\partial\mathbb{D}$ and
\[
z_1+\cdots+z_n=w=0.
\]
For $r=1$, equality holds precisely when $w$ is a zero of $p$ with multiplicity at least $2$.
\end{theorem}

Note that $E$ and $O_n$ are strictly decreasing on $[0,1]$, with
\[
E(0)=O_n(0)=1,\qquad
E(1)=O_n(1)=0,
\]
and $O_n<E$ on $(0,1)$.

We also establish the following elementary geometric proposition, which describes a related extremal configuration.

\begin{proposition}
\label{prop:prop_2}
Let
\[
p(z)=\prod_{j=1}^n(z-z_j),\qquad n\geq 2,
\]
and let $w$ be a critical point of $p$. Then
\begin{equation}
\label{eq:eq_03}
\operatorname{dist}(w,Z(p))
\leq
\frac{1}{n}
\min_{a\in\mathbb{C}}
\sum_{j=1}^n |z_j-a|.
\end{equation}
Equality holds precisely in the following two cases:
\begin{enumerate}
\item[(i)]
\[
|z_1-w|=\cdots=|z_n-w|,
\qquad
\sum_{j=1}^n(z_j-w)=0.
\]

\item[(ii)]
\[
p(z)=(z-\alpha)^k(z-\beta)^{n-k},
\quad 1\le k \le n-1, \quad
w=\frac{(n-k)\alpha+k\beta}{n}.
\]
\end{enumerate}
When $\alpha=\beta$, or when $k=n/2$, case \textup{(ii)} is contained in
case \textup{(i)}.
\end{proposition}

\section{Proof of Proposition~\ref{prop:prop_2}}

If $w\in Z(p)$ and equality holds in~\eqref{eq:eq_03}, then the right-hand
side is zero, so all zeros coincide with $w$. Assume henceforth that
$p(w)\ne0$. Since
\[
\frac{p'(w)}{p(w)}=\sum_{j=1}^n\frac1{w-z_j}=0,
\]
for every $a\in\mathbb C$,
\[
\sum_{j=1}^n\frac{z_j-a}{z_j-w}
=\sum_{j=1}^n\left(1+\frac{w-a}{z_j-w}\right)=n.
\]
Writing $d=\mathrm{dist}(w,Z(p))$, we obtain
\[
n\leq\sum_{j=1}^n\frac{|z_j-a|}{|z_j-w|}
\leq\frac1d\sum_{j=1}^n|z_j-a|,
\]
and minimizing over $a$ proves~\eqref{eq:eq_03}.

Let $a$ be a geometric median and suppose equality holds. Then
\begin{equation}
\label{eq:eq_star}
|z_j-w|=d,\qquad
\frac{z_j-a}{z_j-w} \in [0,\infty)
\quad\text{whenever }z_j\ne a.                 
\end{equation}
If $a\notin Z(p)$, all zeros are equidistant from $w$, and
\[
0=\sum_{j=1}^n\frac1{z_j-w}
=\frac1{d^2}\sum_{j=1}^n \overline{z_j-w},
\]
giving~(i). Conversely, the conditions in~(i) make $w$ both critical and a
geometric median, hence equality holds.

Suppose now $a\in Z(p)$. After translation and rotation, let $w=0$ and
$a=\rho>0$. Since $\rho\geq d$, for every zero $z_j\ne a$, condition~\eqref{eq:eq_star}
gives
\[
|z_j|=d,\qquad 1-\frac{\rho}{z_j}\geq0.
\]
Thus $z_j$ is real. The possibility $z_j=d$ would imply $\rho=d$ and hence
$z_j=a$, so every zero different from $a$ equals $-d$. Thus, unless all
zeros coincide, $p$ has the form in~(ii).

Conversely, if
\[
p(z)=(z-\alpha)^k(z-\beta)^{n-k},\qquad
w=\frac{(n-k)\alpha+k\beta}{n},
\]
then
\[
\operatorname{dist}(w,Z(p))
=\frac{\min\{k,n-k\}}{n}|\alpha-\beta|,
\qquad
\min_{a\in\mathbb C}\sum_{j=1}^n|z_j-a|
=\min\{k,n-k\}|\alpha-\beta|,
\]
so equality holds.

\section{Proof of Theorem 1}

Let $d:=\mathrm{dist}(w, Z(p))$. After a rotation, we may assume that $w=r\in[0,1]$. If $p(r)=0$, the assertion is
immediate; hence assume $p(r)\neq 0$. We have
\[
\frac{p'(r)}{p(r)}
=
\sum_{j=1}^n \frac{1}{r-z_j}
=
0.
\]
If $r=0$, then $\operatorname{dist}(0,Z(p))\le 1$. Equality is attained precisely when all roots $z_i, 1\le i\le n$, of $p$ lie on $|z|=1$ and $p'(0)=0$. This implies 
\[
0=\sum_{j=1}^n \frac{1}{z_j}=\sum_{j=1}^n \overline{z_j},
\]
thus $\sum_{j=1}^n z_j=0$.

Suppose $r=1$. Since $w$ is a critical point and all zeros lie in
$\mathbb D$, the Gauss--Lucas theorem (\cite{Marden}) implies that a critical point on
$\partial\mathbb D$ must be a zero of $p$. Hence
\[
\operatorname{dist}(w,Z(p))=0,
\]
which means
\[
E(1)=O_n(1)=0.
\]
Equality therefore holds precisely when $w$ is a zero of $p$; together with
$p'(w)=0$, this means that $w$ has multiplicity at least $2$.

Assume now $0<r<1$. Let
\[
x_j:=\frac{1}{z_j-r}.
\]
Then
\begin{equation}
\label{eq:eq_1}
\sum_{j=1}^n x_j=0.
\end{equation}

Note that 
\begin{equation}
\label{eq:eq_2}
x_j\in X:= \left\{z\in\mathbb{C}: \left|z-\frac{r}{1-r^2}\right|\ge \frac{1}{1-r^2}\,,\, |z|\le \frac{1}{d} \right\},\, 1\le j\le n.
\end{equation}
Moreover $|z_j|\le 1, 1\le j\le n$ is equivalent to $x_j\in X, 1\le j\le n$. Let us define 
\[
S_n(X):=\{z_1+\cdots+z_n:z_i\in X,\ 1\leq i\leq n\}.
\]
Thus,
\[
R_n(r)=\sup \,\{d: d>0, 0\in S_n(X)\}.
\]
We put 
\[
a=\frac{r}{1-r^2}\,,\quad b=\frac{1}{1-r^2}\,, \quad R=\frac{1}{d}.
\] 
Let $u_{1,2}=\alpha\pm i\beta$ be the two intersection points of the circles $|z-a|=b$ and $|z|=R$. Since $\operatorname{Re} S_n(X) \le n\alpha$, from
\begin{align*}
\sup \,\{d: d>0, 0\in S_n(X)\}&=\left(\inf \,\{R: R>0, 0\in S_n(X)\}\right )^{-1}\\
&\le \left( \inf \, \{R: R>0, \alpha(R)\ge 0\}\right)^{-1},
\end{align*}
we get 
$$R_n(r)\le \sqrt{1-r^2}.$$
In the case $n=2m, m\ge 1$ and $\displaystyle R=\left({1-r^2}\right)^{-1/2}$, taking each of $u_1$ and $u_2$ with multiplicity $m$ yields $0\in S_n(X)$, so $R_n(r)= \sqrt{1-r^2}$.

\smallskip
Assume that $n=2m+1$. Set $d_0:= O_n(r)$ and $R=1/d_0$. For the points $z_j=u_1, 1\le j\le m$, $z_j=u_2, m+1\le j\le 2m $ and $z_{2m+1}=a-b$, we have
\[
m u_1+m u_2+(a-b)=m(u_1+u_2)+a-b=2m\alpha+a-b=0.
\]
Hence $0\in S_n(X)$. Therefore
\[
R_{2m+1}(r)=\sup \,\left\{d: d\ge d_0, 0\in S_n(X)\right\}=\left(\inf \, \left\{R: 0<R\le \frac{1}{d_0}, 0\in S_n(X) \right\}\right)^{-1}.
\]
It suffices to consider only those $R$ for which $\alpha=\alpha(R)>0$, since otherwise $0\notin S_n(X)$. Let
\[
\lambda:=\frac{\alpha+b-a}{\beta}.
\]
A straightforward calculation shows that 
\begin{equation}
\label{eq:eq_5}
\text{if}\quad  R\le \frac{1}{d_0} \quad \text{then} \quad  \lambda\le \frac{\beta}{\alpha}.
\end{equation}
We enlarge $X$ by replacing the two arcs of $|z-a|=b$ joining the points
$a-b$ to $u_1,u_2$ by the corresponding chords. Denote the resulting
symmetric set by $\widetilde X$, thus $X\subset\widetilde X$. 

For
$0\le y\le\beta$ the rightmost point of $\widetilde X$ at height $y$
has real coordinate $a-b+\lambda y$, whereas for $y\ge\beta$ it has real coordinate $\sqrt{R^2-y^2}$.
Hence, the rightmost boundary of the upper half of $\widetilde X$ is given
by $x=F(y)$, where
\[
F(y):=
\begin{cases}
a-b+\lambda y, & 0\le y\le\beta,\\[1mm]
\sqrt{R^2-y^2}, & \beta\le y\le R,
\end{cases}
\]
Suppose that $\displaystyle z_j\in \widetilde{X}, \sum_{j=1}^n z_j\in \mathbb{R}$. Let $p$ be the number of points $z_j, 1\le j\le n$ that lie in $\{z\in \mathbb{C}: \operatorname{Im} z > 0\}$ and let $q=n-p$. Let
\[
T=\sum_{\operatorname{Im}z_j > 0}\operatorname{Im}z_j
 =-\sum_{\operatorname{Im}z_j<0}\operatorname{Im}z_j.
\]
Since $\widetilde X \cap \{z\in \mathbb{C}: \operatorname{Im} z \ge 0\}$ is convex, the maximal real contribution of the $p$ upper points is $\displaystyle pF(T/p)$
and similarly for the $q$ points with $\operatorname{Im} z_j\le 0$. Therefore, for fixed
$p,q,T$,
\[
\operatorname{Re}\sum_{j=1}^n z_j
\le
pF(T/p)+qF(T/q).
\]

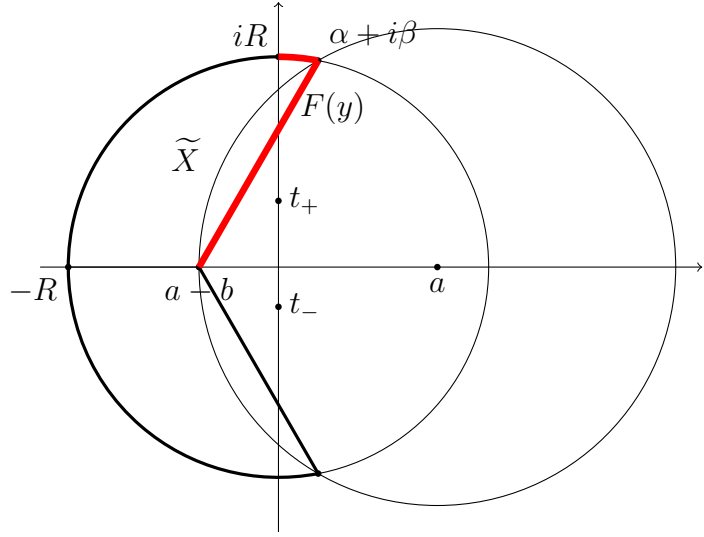
\begin{center}

\begin{tikzpicture}[scale=1.75]

%
%

\def\a{1.2}
\def\b{1.8}
\def\R{1.587451}

\def\al{0.3}
\def\be{1.558846}

\def\tplus{0.5}
\def\tminus{0.3}

\pgfmathsetmacro{\amb}{\a-\b}
\pgfmathsetmacro{\minusR}{-\R}

\pgfmathsetmacro{\theta}{atan2(\be,\al)}


\draw[->] (-1.8,0) -- (3.2,0);
\draw[->] (0,-2.0) -- (0,2.0);


\draw[thin] (0,0) circle[radius=\R];

\draw[thin] (\a,0) circle[radius=\b];


\fill (\a,0) circle (0.7pt);
\fill (\amb,0) circle (0.7pt);
\fill (\minusR,0) circle (0.7pt);

\fill (\al,\be) circle (0.7pt);
\fill (\al,-\be) circle (0.7pt);

\fill (0,\R) circle (0.7pt);


\node[below] at (\a,0) {$a$};
\node[below] at (\amb,0) {$a-b$};
\node[below left] at (\minusR,0) {$-R$};

\node[above right] at (\al,\be)
    {$\alpha+i\beta$};

\node[above left] at (0,\R) {$iR$};

\node at (-0.70,0.85) {$\widetilde{X}$};


\fill (0,\tplus) circle (0.7pt);
\node[right] at (0,\tplus) {$t_+$};

\fill (0,-\tminus) circle (0.7pt);
\node[right] at (0,-\tminus) {$t_-$};


\draw[thin]
    (\minusR,0) -- (\amb,0);


\draw[very thick]
    (\amb,0) -- (\al,\be);

\draw[very thick]
    ({\R*cos(\theta)},{\R*sin(\theta)})
    arc[
        start angle=\theta,
        end angle=360-\theta,
        radius=\R
    ];

\draw[very thick]
    (\al,-\be) -- (\amb,0);


\draw[line width=2.5pt, red]
    (\amb,0) -- (\al,\be);

\draw[line width=2.5pt, red]
    ({\R*cos(\theta)},{\R*sin(\theta)})
    arc[
        start angle=\theta,
        end angle=90,
        radius=\R
    ];

\node[right] at (0.08,1.2) {$F(y)$};

\end{tikzpicture}
\end{center}

Without loss of generality, let $p\le q$ (by reflection across the real axis, if needed). Let
$t_+$ and $t_-$ be the average magnitudes of the imaginary parts of
the upper and lower points, respectively. Since their imaginary parts
cancel,
\[
pt_+=qt_-.
\]
Writing $t=t_+$, we obtain
\[
\operatorname{Re}\sum_{j=1}^n z_j
\le H_p(t):=
pF(t)+qF\left(\frac{p}{q}t\right).
\]
We will show that, for each fixed $p$, $0\le p\le m$, the function
$H_p(t)$, $t\in[0,R]$, attains its maximum at $t=\beta$. We will then optimize over $0\le p\le m$. We have
\[
F'(y)=
\begin{cases}
\lambda, & 0<y<\beta,\\[2mm]
-\dfrac{y}{\sqrt{R^2-y^2}}, & \beta<y<R.
\end{cases}
\]
Since $y/\sqrt{R^2-y^2}$ is increasing and
$\sqrt{R^2-\beta^2}=\alpha$, condition (7) gives
\[
F'(y)\le -\frac{\beta}{\alpha}\le -\lambda,
\qquad \beta<y\le R.
\]
Hence, for $t<\beta$, since $p\le q$,
both $t$ and $(p/q)t$ lie in $[0,\beta]$, and
\[
H_p'(t)
=pF'(t)+pF'\left(\frac pq t\right)
=2p\lambda \ge 0.
\]
On the other hand, for $t>\beta$,
\[
F'(t)\le-\lambda,
\qquad
F'\left(\frac pq t\right)\le\lambda,
\]
and therefore
\[
H_p'(t)
=
p\left(
F'(t)+F'\left(\frac pq t\right)
\right)\le0.
\]
Thus, $H_p(t)$ is increasing on $[0,\beta]$ and is decreasing on $(\beta, R]$. Hence its maximum is attained at $t=\beta$, and
\[
H_p(\beta)
=2p\alpha+(q-p)(a-b)
=n(a-b)+2p(\alpha+b-a).
\]
Since $\alpha+b-a>0$, it follows that $H_p(\beta)$ is increasing in $p$. Taking into account $p\le m<q$, the maximum is attained for $p=m$, hence $q=m+1$. Consequently
\begin{equation}
\label{eq:eq_6}
\operatorname{Re}\sum_{j=1}^n z_j \le H_m(\beta) = 2m\alpha(R)+(a-b)= (n-1)\alpha(R)+(a-b).
\end{equation}
Therefore $0\in S_n(X)$ requires $2m\alpha(R)+(a-b)= (n-1)\alpha(R)+(a-b)\ge 0$, and when $(n-1)\alpha(R)+(a-b) = 0$ we have $0\in S_n(X)$. Since $(n-1)\alpha(R)+(a-b)$ is an increasing function of $R$, we get  
\[
\inf \left\{R: 0<R\le \frac{1}{d_0}, 0\in S_n(X) \right\}= \inf \left\{R: 0<R\le \frac{1}{d_0}, (n-1)\alpha+(a-b)=0\right\}
\]
Therefore $R_n(r)=O_n(r)$. For $z_j\in \widetilde{X}$, equality in \eqref{eq:eq_6} occurs only when the $m$ points in $\operatorname{Im} z>0$ are equal to $u_1=\alpha+i \beta$, while the remaining points lie on the chord $[a-b, u_2]$, with $t_-=m\beta/(m+1)$. 

Since $X\subset\widetilde{X}$ meets the chord $[a-b,u_2]$ only at its endpoints, equality for $z_j\in X$ occurs precisely when the configuration consists of $m$ copies of $u_1$, $m$ copies of $u_2$, and the point $a-b$.

Transforming back by $z=r+1/x$, the point $a-b$ gives the zero $-1$, while $u_1,u_2$ give the zeros $e^{\pm i\theta}$, with $\cos\theta$ as stated in Theorem~\ref{thm:thm_1}.


\bigskip

\begin{thebibliography}{9}

\bibitem{Marden}
M.~Marden,
\textit{Geometry of Polynomials},
2nd ed., Math. Surveys, No.~3,
Amer. Math. Soc., Providence, RI, 1966.

\bibitem{Mazur}
L.~Mazur,
\textit{A Computer-Assisted Proof of Sendov's Conjecture},
public version of 5 August 2026, ProofAtlas, 2026.

\bibitem{Tao}
T.~Tao,
\textit{A digestion of the proof of Sendov's conjecture},
What's new, 12 August 2026.

\bibitem{SofiShah}
G.~M. Sofi and W.~M. Shah,
\textit{Distribution of zeros and critical points of a polynomial, and Sendov's conjecture},
J. Contemp. Math. Anal. \textbf{58} (2023), 384--388.

\end{thebibliography}
\end{document}